\documentclass{amsart} 
\usepackage{amsmath,amssymb,amsthm}
\input diagrams
\numberwithin{equation}{section}
\newtheorem{theorem}[equation]{Theorem}
\newtheorem{proposition}[equation]{Proposition}

\newtheorem{lemma}[equation]{Lemma}

\theoremstyle{remark}

\newtheorem{remark}[equation]{Remark}
\newtheorem*{note}{Note}
\newtheorem*{notation}{Notation}
\newtheorem*{remark*}{Remark}
\begin{document}
\title[Vector bundles on non-K\"{a}hler elliptic surfaces]
{Holomorphic rank-2 vector bundles on non-K\"{a}hler elliptic surfaces}

\author{Vasile Br\^{\i}nz\u{a}nescu}
\address{Institute of Mathematics "Simion Stoilow",
Romanian Academy, P.O.Box 1-764, RO-70700,
Bucharest, Romania}
\email{Vasile.Brinzanescu@imar.ro} 

\thanks{The first author was partially supported by grant CNCSIS 33518/2002
and by EURROMMAT Programme FP5/2000-2003}

\author{Ruxandra Moraru}
\address{Department of Mathematics and Statistics, Burnside Hall,
McGill University, 805 Sherbrooke Street West,
Montreal, Quebec, Canada, H3A 2K6}
\email{moraru@math.mcgill.ca}

\dedicatory{\it Dedicated to Professor F. Hirzebruch for his 75th birthday.}


\begin{abstract}
The existence problem for vector bundles on a smooth compact complex surface 
consists in determining which 
topological complex vector bundles admit holomorphic structures. For 
projective surfaces, Schwarzenberger proved that a topological complex 
vector bundle admits a holomorphic (algebraic) structure if and only if its 
first Chern class belongs to the Neron-Severi group of the surface. In contrast, for 
non-projective surfaces there is only a necessary condition for the
existence problem (the discriminant of the vector bundles must be positive)
and the difficulty of the problem resides in the lack of a 
general  method for constructing non-filtrable vector bundles. In this paper, 
we close the existence problem in the rank-2 case, by giving necessary and 
sufficient conditions for the existence of holomorphic rank-2 vector 
bundles on non-K\" ahler elliptic surfaces.  
\end{abstract}
\maketitle

\section{Introduction}

In this paper, we study the existence of holomorphic vector bundles on 
non-K\"{a}hler elliptic surfaces; their classification and
stability properties are discussed in 
\cite{Brinzanescu-Moraru1,Brinzanescu-Moraru2}. 
Let $X$ be a smooth compact complex surface.
The existence problem for vector bundles on $X$
consists in determining which 
topological complex vector bundles admit holomorphic structures, or 
equivalently, in finding all triples $(r,c_1,c_2)$ in $\mathbb{N}\times NS(X) 
\times \mathbb{Z}$ for which there exists a rank-$r$ holomorphic 
vector bundle on $X$ with Chern classes $c_1$ and $c_2$. 
For projective surfaces, Schwarzenberger \cite{Schw} proved that any triple 
$(r,c_1,c_2)$ in $\mathbb{N} \times NS(X) \times \mathbb{Z}$
comes from a rank-$r$ holomorphic (algebraic) vector bundle. In contrast, 
for non-projective surfaces, there is a natural necessary condition for 
the existence problem \cite{BaL,BrF,LeP}:
\[ \Delta (r,c_1,c_2):= \frac{1}{r} \left( c_2 - \frac{r-1}{2r}c_1^2 
\right) \geq 0. \]
One can always construct filtrable bundles by using extensions of coherent sheaves; 
in fact, on a non-algebraic surface $X$, 
there exists a filtrable rank-$r$ holomorphic vector bundle $E$ with Chern classes $c_1$ 
and $c_2$ if and only if its
discriminant $\Delta (E)$ satisfies the inequality
\[ \Delta (E):= \Delta (r,c_1,c_2) \geq m(r,c_1), \]
where
\[ m(r,c_1):= -\frac{1}{2r}\mbox{max} \left\{ \sum_1^r \left( \frac{c_1}{r}-
\mu_i \right)^2,\mu_1,\dots,\mu_r\in NS(X), \sum_1^r\mu_i=c_1 \right\} \]
(see \cite{BaL,BrF,LeP}).
Therefore, the only unknown situations occur for bundles of rank greater than one
that have a discriminant in the interval $\left[ 0 , m(r,c_1) \right)$;
vector bundles with such discriminants will, of course, be non-filtrable and the difficulty 
of the problem resides in the lack of a general method for constructing non-filtrable bundles. 
One is thus compelled to focus on particular classes of surfaces, 
to find specific construction methods.

The existence of bundles on non-projective surfaces
is, in general, still an open question, which has been completely settled only in the case of 
primary Kodaira surfaces \cite{ABT}. 
For rank-$2$ holomorphic vector bundles, the problem has been solved for 
complex 2-tori \cite{To}, as well as for surfaces of class VII and K3 surfaces 
\cite{TT}; since the method used in \cite{TT} 
(Donaldson polynomials) seems to also work for (non-algebraic) 
K\" ahler elliptic surfaces, only the case of general non-K\" ahler elliptic surfaces 
remains.
In this article, we close 
the existence problem in the rank-$2$ case, by giving necessary and 
sufficient conditions for the existence of holomorphic rank-$2$ vector 
bundles on non-K\" ahler elliptic surfaces.

Recall that a surface is said to be elliptic if it admits a holomorphic 
fibration over a curve with generic fibre an elliptic curve;
for instance, non-K\" ahler elliptic surfaces are given by holomorphic fibrations
without a section whose smooth fibres are isomorphic to a fixed elliptic curve.
For vector bundles on any elliptic fibration $\pi: X \rightarrow B$, 
restriction to a fibre is a natural operation:
there exists a divisor in the relative Jacobian 
$J(X)$ of $X$, 
called the {\em spectral curve} or {\em cover} of the bundle, 
that encodes the isomorphism
class of the bundle over each fibre of $\pi$. 
This divisor is an important invariant of bundles on elliptic fibrations,
which has proven very useful in their study  
(see \cite{F1,FM,FMW,BJPS,D}) for projective fibrations, \cite{DOPW1,DOPW2}
for Calabi-Yau threefolds without a section, and \cite{B-H,Moraru,T}
for non-K\"{a}hler fibre bundles).
The spectral construction presented in this paper is a modification
of the Fourier-Mukai transform for 
certain elliptic fibrations without a section,
which will be used in \cite{Brinzanescu-Moraru1} to define
a twisted Fourier-Mukai transform that is specific to non-K\"{a}hler elliptic surfaces. 

The paper is organised as follows. 
We begin by presenting and proving some topological and geometrical 
properties of non-K\" ahler elliptic surfaces;
in particular, we show that if $\pi: X \rightarrow B$ is such a surface, 
then the restriction of {\em any} vector bundle on $X$ to a smooth fibre of $\pi$ {\em always} 
has degree zero.
Unlike the algebraic case \cite{FM},
the description of line bundles on non-K\" ahler elliptic surfaces is not straightforward;
indeed, even though these surfaces have very few divisors
(they are given by the fibres of $\pi$), 
there exist many line bundles on them. 
Nonetheless, we are able to establish a correspondence between line 
bundles on a non-K\" ahler elliptic surface and sections of its 
relative Jacobian; this follows from results of 
\cite{B1,B2,B3,BrU} regarding the Neron-Severi and Picard groups of these surfaces. 
In the third section, we extend the spectral construction of \cite{B-H,Moraru}
to the case of holomorphic rank-$2$ vector bundles on non-K\"{a}hler 
elliptic surfaces.
Finally, the last section contains the proof of the existence 
theorems.  

We end by noting that the techniques developed here and in 
\cite{Brinzanescu-Moraru1,Brinzanescu-Moraru2} can be used to
solve existence and classification problems for holomorphic vector bundles 
of arbitrary rank on non-K\"{a}hler elliptic and torus fibrations.  

\section{Line bundles}
\label{line bundles}

Let $X\stackrel{\pi}{\rightarrow}B$ be a minimal non-K\"ahler 
elliptic surface, with $B$ a smooth compact connected curve; it is 
well-known that $X \stackrel{\pi}{\rightarrow} B$ is a quasi-bundle over $B$, 
that is, all the smooth fibres of $\pi$ are isomorphic to a fixed elliptic curve $T$ 
and the singular ones (if any) are isogeneous to multiples of $T$ (see \cite{Kod,B3}). We begin by 
presenting several topological and geometric properties of these surfaces.
 
Let $T^*$ denote the dual of $T$ 
(we fix a non-canonical identification $T^*:= \mbox{Pic}^0(T)$). In this 
case, the Jacobian surface associated to $X\stackrel{\pi}{\rightarrow}B$ 
is simply
\[ J(X)=B\times T^*\stackrel{p_1}{\rightarrow}B \]
(see, for example, \cite{Kod,BPV,B1}) and the surface 
is obtained from its relative Jacobian by a finite number of logarithmic transformations 
\cite{Kod,BPV,BrU}. 
Also, if $X$ has multiple fibres $T_1,\dots, T_r$, with corresponding 
multiplicities $m_1,\dots, m_r$, then its canonical bundle is given by
\[K_X = \pi^\ast K_B \otimes \mathcal{O}_X\left(\sum_{i=1}^r (m_i-1)T_i\right). \] 
Finally, we the following identification \cite{B1,B2,BrU}:
\[NS(X)/Tors(NS(X)) \cong Hom(J_B, Pic^0(T)),\]
where $NS(X)$ is the the Neron-Severi group of the surface and
$J_B$ denotes the Jacobian variety of $B$; the torsion of 
$H^2(X,\mathbb{Z})$ is generated by the classes of the fibres (both smooth 
and multiple). In the remainder, the class modulo $Tors(H^2(X,\mathbb{Z}))$ 
of an element $c \in H^2(X,\mathbb{Z})$ will be denoted $\widehat{c}$. 
Given these considerations, we have:

\begin{lemma}\label{c^2}
Let $X \stackrel{\pi}{\rightarrow} B$ be a non-K\"{a}hler elliptic surface.

(i) If $c \in NS(X)$, then $\pi_\ast(c) = 0$ and $c \cdot \beta = 0$ 
for all $\beta \in Tors(H^2(X,\mathbb{Z}))$.

(ii) For any element $c \in NS(X)$, $c^2 = - 2\deg(\widehat{c})$.
\end{lemma}

\begin{proof} 
The lemma is certainly true for torsion classes. Let us then 
assume that $c \notin Tors(NS(X))$ and choose a line bundle $L$ on $X$ 
with first Chern class $c$. Then $\widehat{c} \neq 0$ and, by fixing a 
base-point in $B$, the cohomology class $\widehat{c}$ can be considered as 
a covering map $\widehat{c}: B \rightarrow Pic^0(T)$ such that
\[ {\widehat{c}}^{\ -1}(\lambda_0) = \{ b \in B \ | \ L|_{F_b} 
\simeq \lambda_0 \}. \]
Since $\widehat{c} \neq 0$, we have 
${\widehat{c}}^{\ -1}(\mathcal{O}_T) \neq B$. Therefore, the stalk of 
$\pi_\ast L $ is zero at the generic point in $B$ and the direct image sheaf $\pi_\ast L$ vanishes;
furthermore, the higher direct image sheaf $R^1\pi_\ast L $ is a torsion sheaf supported on 
${\widehat{c}}^{\ -1}(\mathcal{O}_T)$. In particular, 
$\pi!L = -R^1 \pi_\ast L$ and, by Grothendieck-Riemann-Roch, the pushdown
$\pi_\ast(c)$ is equal to the rank of the torsion sheaf $R^1\pi_\ast L$, 
which is zero. Let $\beta$ be a generator of the torsion of 
$H^2(X,\mathbb{Z})$. The class of $\beta$ is then the first Chern class 
of a sheaf on $X$ of the form $\mathcal{O}(F)$, where $F$ is 
a fibre of $\pi$ of multiplicity $m \geq 1$. Consequently, 
the pullback to $X$ of the positive generator $h$ of $H^2(B,\mathbb{Z})$
is equal to $c_1(\mathcal{O}(mF))$ and, by the Projection formula, we have
\[ m(c \cdot \beta)h = \pi_\ast(c \cdot \pi^\ast h) 
= \pi_\ast(c) \cdot h = 0,\]
that is, $c \cdot \beta = 0$, proving (i). Combining the 
results of (i) with Grothendieck-Riemann-Roch, we obtain 
$c_1(R^1\pi_\ast L) = - \frac{1}{2}c^2 \cdot h$. Hence, the degree of the
map $\widehat{c}$ is equal to $\# ({\widehat{c}}^{\ -1}(\mathcal{O}_T)) = -\frac{1}{2}c^2$
and we are done.
\end{proof}

\begin{lemma}\label{degree of restriction of line bundle}
Let $\pi :X\rightarrow B$ be a non-K\"{a}hler elliptic surface and 
$\mathcal{L}$ a line bundle on $X$. The restriction of $\mathcal{L}$ 
to any smooth fibre 
of $\pi$ has degree zero.
\end{lemma}

\begin{proof}
Let $m_1T_1,\;m_2T_2,...,\;m_{\ell}T_{\ell}$ be the
multiple fibres of $\pi$ and set $b_i=\pi (T_i)$. Denote $m$ the least 
common multiple of $m_1,m_2,\dots,m_{\ell}$ and choose a non-negative integer 
$e$ such that $m$ divides $\ell +e$; next, take distinct points 
$b_{\ell +1},\dots, b_{\ell +e}$, which are different from $b_i$, $i=1,\dots,\ell$,
and fix a point $b$ with $T_b$ smooth. Then, there exists at least one 
line bundle $M$ on $B$ with the property that
\[M^{\otimes m}\cong \mathcal{O}_B(b_1 + ... +b_{\ell +e});\]
such a line bundle defines an $m$-cyclic covering $\varepsilon :B'\rightarrow B$
that is totally ramified at $b_1,...,\;b_{\ell +e}$ (see \cite{BPV},  
Chapter I, Lemma 17.1). By Lemma 3.18 in \cite{B3}, there 
exists a principal $T$-bundle $\pi ':X'\rightarrow B'$ and an $m$-cyclic covering 
$\psi :X'\rightarrow X$ over 
$\varepsilon :B'\rightarrow B$; let $\tilde T$ be a connected component of 
$\psi^{-1}(T_b)$. Then $\tilde T$ is a fibre of $\pi '$ and the restriction 
$\tilde T\rightarrow T_b$ of $\psi$ is an isomorphism. Therefore, we have
\[c_1(\mathcal{L}|_{T_b})=c_1(\psi^*(\mathcal{L})|_{\tilde T})=0,\]
because $\pi ':X'\rightarrow B'$ is a principal elliptic bundle \cite{B3,T}.
\end{proof}

\begin{remark*}
Similar results are stated in \cite{ABT,T} for non-K\"{a}hler principal elliptic bundles,
that is, non-K\"{a}hler elliptic surfaces without multiple fibres.
\end{remark*}

Referring to Lemma \ref{degree of restriction of line bundle} and \cite{BrU}, we can therefore
associate to any line bundle $\mathcal{L}$ on $X$ a holomorphic mapping 
$\varphi :B\rightarrow T^*$ such that
\[\mathcal{L}|_{T_b}=\varphi (b),\]
for any smooth fibre $T_b$, that is, a section of $J(X) = B \times T^\ast$.
Conversely, one can associate to every section of $J(X)$ a line bundle on $X$, as stated in:
\begin{proposition}\label{section of line bundle}
Let $\pi :X\rightarrow B$ be a non-K\"{a}hler elliptic surface, with
general fibre $T$, and $J(X)=B\times T^*$ be the associated Jacobian surface 
of $X$. Then:

(i) For any section $\Sigma \subset J(X)$, there exists a line bundle 
$\mathcal{L}$ on $X$ whose restriction to every smooth fibre $T_b$ is the same 
as the line bundle $\Sigma_b$ of degree zero on $T=T_b$.

(ii) The set of all line bundles on $X$ that restrict, on every smooth 
fibre of $\pi$,  
to the line bundle of degree zero determined by the section 
$\Sigma$ is a principal homogeneous space over $P_2$, where $P_2$ is the 
subgroup of line bundles on $X$ generated by $\pi^*\mbox{Pic}(B)$ and the
$\mathcal{O}_X(T_i)$'. 
\end{proposition}

\begin{proof} 
Choose a general point $b\in B$ with $T_b$ smooth and consider the 
natural restriction morphism
$r:\rm{Pic}(X)\rightarrow \rm{Pic}(\pi^{-1}(b))=\rm{Pic}(T).$
Let $(P_j)$ be the filtration of $\mbox{Pic}(X)$ defined by
\[ \mbox{$P_0=\rm{Pic}(X)$, $P_1=\rm{Ker}(r)$, and $P_2$}. \]
Set $N(X):=P_0/P_1$ and $\tilde{N}(X):=\{c_1(L)|\;L\in N(X)\}$. 
Referring to \cite{B1} and \cite{BrU}, we have ${\tilde N}(X)=0$ and
\[NS(X)/\mbox{Tors}NS(X)\cong \mbox{Hom}(J_B,T^*)\cong P_1/P_2.\]
Consequently, $N(X)\subset \mbox{Pic}^0(T)$. Since any line bundle in 
$\mbox{Pic}^0(T)$ is invariant by translations, we obtain
\[N(X)=\mbox{Pic}^0(T)\]
by Lemma 
\ref{degree of restriction of line bundle} and \cite{BrU}.
Let $\lambda =\Sigma_b\in T^*$ 
and let $\Sigma^\lambda$ be the constant section 
$B\times \lambda \subset J(X)$. Following the construction in \cite{BrU}, 
the line bundle $\lambda \in T^*$ extends to a line bundle $\mathcal{L}^\lambda$ 
on $X$ that corresponds to the 
constant section $\Sigma^\lambda$. Let $B_0$ be the zero section of $J(X)$. 
Given the identification $P_1/P_2\cong \mbox{Hom}(J_B,T^*)$,
there exists a line bundle $\mathcal{L}_1$ in $P_1=\mbox{Ker}(r)$ 
whose corresponding element in $\mbox{Hom}(J_B,T^*)$ is a section that is linearly 
equivalent to $\Sigma - \Sigma^\lambda + B_0$ (look at the addition law 
of the group $\mbox{Hom}(J_B,T^*)$). The line bundle 
$\mathcal{L}=\mathcal{L}_1\otimes \mathcal{L}^\lambda$ is then such that its 
restriction to every smooth fibre $T_b$ is the same as the line bundle 
$\Sigma_b\in T^*$, proving (i). If the line bundles $\mathcal{L}'$ and $\mathcal{L}$ on 
$X$ both have the above property, then by the same isomorphism, 
$\mathcal{L}'\otimes \mathcal{L}^{-1}\in P_2$ and we are done.
\end{proof}

We can now characterise the sections of the Jacobian surface as follows.
\begin{lemma}\label{degree of section} 
Let $X$ be a non-K\"{a}hler elliptic surface. 
Then, any section $\Sigma$ of the Jacobian surface $J(X)$ of $X$
has trivial self-intersection. Furthermore, if $\mathcal{L}$ is any line bundle on $X$ corresponding 
to the section $\Sigma$ of $J(X)$, then
\[\Sigma  \cdot B_0 = -c_1^2(\mathcal{L})/2,\]
where $B_0$ denotes the zero section of $J(X)$.
\end{lemma}

\begin{proof}
The invariants of the Jacobian surface $J(X)=B\times T^*$ are
\[ \mbox{$p_g(J(X))=g$, $q(J(X))=g+1$, and $K_{J(X)}=p_1^*K_B$},\]
where $g$ is the genus of the curve $B$; the adjunction formula gives $\Sigma^2=0$.
Let $\hat{c}_1$ be the class of $c_1(\mathcal{L})$ in $NS(X)/Tors(NS(X)) \cong Hom(J_B,T^\ast)$.
As in the proof of Lemma \ref{c^2}, 
we can then think of $\hat{c}_1$ as being a covering map 
$\hat{c}_1:B\rightarrow T^*$ of degree $-c_1^2(\mathcal{L})/2$;
since the degree of $\hat{c}_1$ is also equal to $\Sigma \cdot B_0$, the lemma follows. 
\end{proof}

We end the section by giving a description of torsion line bundles on a principal elliptic 
bundle $X \stackrel{\pi}{\rightarrow} B$; the surface is now isomorphic to a quotient of the form
\[ X = \Theta^\ast / \langle \tau \rangle ,\]
where $\Theta$ is a line bundle on $B$ with positive Chern class $d$, 
$\Theta^\ast$ is the complement of the zero section in the total space of 
$\Theta$, and $\langle \tau \rangle$ is the multiplicative cyclic group 
generated by a fixed complex number $\tau$, with 
$| \tau |$ greater than $1$. The standard fibre of this bundle is
\[ T \cong \mathbb{C}^\ast/\langle \tau \rangle \cong \mathbb{C}/
(2\pi i\mathbb{Z} + \ln(\tau)\mathbb{Z}).\]
(We assume $d$ to be positive so that the surface $X$ is non-K\"{a}hlerian.) 

The set of all holomorphic line bundles on $X$ with trivial Chern class 
is given by the zero component of the Picard group $\text{Pic}^0(X)$. 
Referring to Proposition 1.6 in \cite{T}, one has
\[ \text{Pic}^0(X) \cong \text{Pic}^0(B) \times \mathbb{C}^\ast. \]
Any line bundle in $\text{Pic}^0(X)$ is therefore of the form 
$H \otimes L_\alpha$, where $H$ is the pullback to $X$ of an element of 
$\text{Pic}^0(B)$ and $L_\alpha$ is the line bundle corresponding to the 
constant automorphy factor $\alpha \in \mathbb{C}^\ast$. We illustrate 
this by constructing the restriction of the universal (Poincar\'{e}) 
line bundle $\mathcal{U}$ over $X \times \text{Pic}^0(X)$ to 
\[X \times \mathbb{C}^\ast := X \times \{ 0 \} \times \mathbb{C}^\ast .\] 
One starts with a trivial line bundle $\bar{\mathbb{C}}$ on 
$\Theta^\ast \times \mathbb{C}^\ast$ and applies to it the following 
$\mathbb{Z}$-action
\[ \begin{array}{rcl}
\Theta^\ast \times \mathbb{C}^\ast \times \mathbb{Z} & \rTo & 
\Theta^\ast \times \mathbb{C}^\ast \\
(z,\alpha,n) & \rMapsto & (\tau^nz,\alpha).
\end{array} \]
Since this action is trivial on $\mathbb{C}^\ast$, the Poincar\'{e} 
line bundle $\mathcal{U}$ is obtained by identifying 
$s \in \bar{\mathbb{C}}_{(z,\alpha)}$ 
with $\alpha s \in \bar{\mathbb{C}}_{(\tau z,\alpha)}$.

\begin{notation} 
In the remainder, we shall denote by 
$L_\alpha$ the line bundle corresponding to the automorphy factor 
$\alpha \in \mathbb{C}^\ast$.
\end{notation}

\begin{remark*} 
Although the line bundle $L_{\tau^m}$ is trivial over the fibres of $\pi$, one cannot 
define an action of $\mathbb{Z}$ on $\mathbb{C}^\ast$ that leaves the
restriction of the Poincar\'{e} line bundle $\mathcal{U}$ to 
$X \times \mathbb{C}^\ast$ 
invariant. Indeed, if $\mathbb{Z}$ acts on $\mathbb{C}^\ast$, then 
multiplication by $\tau$ is defined on the fibres of $\bar{\mathbb{C}}$ by
\begin{equation}\label{action2}
\begin{array}{rcl}
\tau : \Theta^\ast \times \mathbb{C}^\ast \times \mathbb{C} & \rTo & 
\Theta^\ast \times \mathbb{C}^\ast \times \mathbb{C}\\
(z,\alpha,t) & \rMapsto & (\tau z,\tau\alpha,\alpha t).
\end{array}
\end{equation}
On the surface $X$, $z$ and $\tau z$ define the same point $x$. 
However, \eqref{action2} indicates that $\tau$ sends 
$\mathcal{U}_{(x,\alpha)}$ to 
$\mathcal{U}_{(x,\tau \alpha)} \otimes L_{\tau^{-1},x}$. Hence, the Poincar\'{e} 
line bundle is not invariant under such an action.
\end{remark*}

\section{Holomorphic vector bundles}
\label{spectral curve}

Consider a pair $(c_1,c_2)$ in $NS(X) \times \mathbb{Z}$. Its 
corresponding {\em discriminant} is then given by
\[ \Delta(2,c_1,c_2):=\frac{1}{2} \left( c_2 - \frac{c_1^2}{4} \right)\geq 0.\]
Let $E$ be a rank 2 vector bundle over $X$, with $c_1(E) = c_1$ and 
$c_2(E) = c_2$. 
We fix the following notation:
\[ \mbox{$\Delta(E) := \Delta(2,c_1,c_2)$ and
$n_E := -ch_2(E)$}, \]
where $ch_2(E) = c_1^2/2 - c_2$ is the second Chern character of $E$. 
\begin{remark}\label{n}
Referring to Lemma \ref{c^2}, if $\Delta(2,c_1,c_2) \geq 0$, then 
$n_E \geq 0$. 
\end{remark}

To study bundles on $X$, one of our main tools will be restriction 
of the bundle to the smooth fibres $\pi^{-1}(b) \cong T$ of the fibration 
$\pi : X \rightarrow B$. Since the restriction of any bundle on $X$ to 
a fibre $T$ has first Chern class zero, we consider $E$ as family of 
degree zero bundles over the elliptic curve $T$, parametrised by $B$. 
Given a rank two bundle over $X$, its restriction to a generic 
fibre of $\pi$ is semistable. More precisely, we have:
\begin{proposition} 
Let E be a rank 2 holomorphic vector bundle over $X$. Then,
$E|_{\pi^{-1}(b)}$ is unstable on at most an isolated set of points 
$b \in B$.
\end{proposition}

\begin{proof} 
Suppose that $b\in B$ is a point such that 
$E|_{\pi^{-1}(b)}$ is unstable, splitting as 
$\lambda_b \oplus (\lambda_b')^\ast$
for some line bundles $\lambda_b$ and $\lambda_b'$ in $\mbox{Pic}^{-k}(T)$, $k>0$. 
Consider the elementary modification
\[ 0\rightarrow E'\rightarrow E\rightarrow j_*\lambda_b\rightarrow 0,\]
where $j:T_b\rightarrow X$ is the natural inclusion. Referring to \cite{F2} 
(Chapter II, Lemma 16), the discriminant of $E'$ is given by
\[ \Delta (E')=\Delta (E) + \frac{1}{2}j_*c_1(\lambda_b); \]
furthermore,
\[ \Delta(E')<\Delta (E)\]
because $\deg(\lambda_b) = -k < 0$. 
Therefore, since the existence of $E'$ implies that its discriminant is a non-negative number, 
the result follows.
\end{proof}

\begin{note} These isolated points are called the {\em jumps} of the bundle $E$.
\end{note}

\subsection{The spectral curve of a rank-2 vector bundle}

Let us assume for a moment that $X$ does not have multiple fibres. 
Choose a line bundle $L$ in $\text{Pic}^0(X)$ such that 
$h^0(\pi^{-1}(b),L^\ast \otimes E)$ is zero, for generic $b$. 
The direct image sheaf 
$R^1\pi_\ast(L^\ast \otimes E)$ is therefore a torsion sheaf supported on 
isolated points $b$ such that $E|_{\pi^{-1}(b)}$ is semistable and has $L|_{\pi^{-1}(b)}$ 
as a subline bundle, or $E|_{\pi^{-1}(b)}$
is unstable; consequently, if $h$ is the positive generator of $H^2(B,\mathbb{Z})$, 
then
\[ c_1(R^1\pi_\ast(L^\ast \otimes E)) = - \pi_\ast(ch(E) \cdot td(X)) 
\cdot td(B)^{-1} = n_Eh.\]
However, since the discriminant of $E$ is a non-negative number, 
then so is the integer $n_E$ (see remark \ref{n}): 
the sheaf $R^1\pi_\ast(L^\ast \otimes E)$ is supported on $n_E$ 
points, counting multiplicity.

To obtain a complete description of the restriction of $E$ to 
the fibres of $\pi$, this construction must be repeated for every 
line bundle on $X$; this is done by taking the direct image 
$R^1\pi_\ast$ for all line bundles simultaneously. Let $\pi$ 
also denote the projection
$\pi := \pi \times id: X \times \text{Pic}^0(B) \times \mathbb{C}^\ast \rightarrow B
\times \text{Pic}^0(B) \times \mathbb{C}^\ast$,
where $id$ is the identity map on $\text{Pic}^0(B) \times \mathbb{C}^\ast$,
and let 
$s:X \times \text{Pic}^0(B) \times \mathbb{C}^\ast \rightarrow X$
be the projection onto the first factor. If $\mathcal{U}$ is the universal 
(Poincar\'{e}) line bundle over 
$X \times \text{Pic}^0(B) \times \mathbb{C}^\ast$, one defines
\[ \widetilde{\mathcal{L}} := R^1\pi_\ast (s^\ast E \otimes \mathcal{U}) . \]
This sheaf is supported on a divisor $\widetilde{S_E}$ that is defined 
with multiplicity. We have the following remarks:

\begin{itemize}
\item
Let $H$ be the pullback to $X$ of a line bundle of degree 
zero on $B$. The restriction of $H$ to any fibre $T$ is then trivial, 
implying that the support of 
\[ R^1\pi_\ast (s^\ast E \otimes \mathcal{U} \otimes H)\]
is also $\widetilde{S_E}$. We can therefore restrict the above construction to 
$X \times \mathbb{C}^\ast := X \times \{ 0 \} \times \mathbb{C}^\ast$. 
In the remainder, we will use the same notation for this restriction.

\item
Consider the $\mathbb{Z}$-action on $B \times \mathbb{C}^\ast$ 
induced from the one on $X \times \mathbb{C}^\ast$. For any 
$(b,\alpha)$ in $B \times \mathbb{C}^\ast$, multiplication by $\tau$ sends 
the stalk $\widetilde{\mathcal{L}}_{(x,\alpha)}$ to 
$\widetilde{\mathcal{L}}_{(x,\tau\alpha)} \otimes L_{\tau^{-1},x}$, 
leaving the support of $\widetilde{\mathcal{L}}$ unchanged. 
\end{itemize}
 
By the above remarks, since the quotient $\mathbb{C}^\ast/\langle \tau \rangle$
of $\mathbb{C}^\ast$ by the $\mathbb{Z}$-action is isomorphic to $T^\ast$, 
the support $\widetilde{S_E}$ of $\widetilde{\mathcal{L}}$ descends to a divisor $S_E$ 
in $J(X) = B \times T^\ast$ of the form
\[ S_E := \left( \sum_{i=1}^k \{ x_i \} \times T^\ast \right) 
+ \overline{C},\]
where $\overline{C}$ is a bisection of $J(X)$ (that is, $S_E.T^*=2$ for 
any fibre $T^*$ of $J(X)$) and $x_1, \cdots, x_k$ are points 
(counted with multiplicities) in $B$ that correspond to the jumps of $E$.

If the fibration $\pi$ has multiple fibres, 
the spectral cover of a bundle $E$ on $X$ is then constructed as follows. 
Referring to the proof of Lemma \ref{degree of restriction of line bundle},
there exists a principal $T$-bundle 
$\pi ':X'\rightarrow B'$ over an $m$-cyclic covering $\varepsilon :B'\rightarrow B$. 
Note that the map $\varepsilon$ induces natural $m$-cyclic coverings $\psi :X'\rightarrow X$ 
and $J(X')\rightarrow J(X)$. 
By replacing $X$ with $X'$ (which has no multiple fibres) in the above construction,
we obtain the a spectral cover 
$S_{\psi^\ast E}$ of $\psi^\ast E$ as a divisor in $J(X')$.
We define the spectral cover $S_E$ of $E$ as the projection of 
$S_{\psi^\ast E}$ in $J(X)$; one easily sees that $S_E$ does indeed give 
the isomorphism type of $E$ over each smooth fibre of $\pi$.

\begin{remark*}
The above construction can be defined for any rank-$r$ vector bundle. 
In particular, for a line
bundle, the spectral cover corresponds to the section of the 
Jacobian surface $J(X)$ defined in 
section 2.
\end{remark*}

\subsection{The graph of a rank-2 vector bundle}
\label{graph}
Let $\delta$ be the determinant line bundle of $E$. It then defines the
following involution on the relative Jacobian 
$J(X) = B \times T^\ast$ of $X$:
\[ \begin{array}{rcl}
i_\delta: J(X) & \rightarrow & J(X) \\
(b,\lambda) & \mapsto & (b,\delta_b \otimes \lambda^{-1}), 
\end{array} \]
where $\delta_b$ denotes the restriction of $\delta$ to the fibre 
$T_b = \pi^{-1}(b)$. For fixed point $b$ in $B$, the involution induced
on the corresponding fibre of $p_1: J(X) \rightarrow B$
has four fixed points (the solutions of $\lambda^2 = \delta_b$). 
Taking the quotient of $J(X)$ by this involution, each fibre
of $p_1$ becomes $T^\ast/i_\delta \cong \mathbb{P}^1$ and the
quotient $J(X)/i_\delta$ is isomorphic to a ruled surface 
$\mathbb{F}_\delta$ over $B$. 
Let $\eta: J(X) \rightarrow \mathbb{F}_\delta$ be the canonical map.
By construction, the spectral curve $S_E$ associated to $E$ is invariant
under the involution $i_\delta$ and descends to the quotient 
$\mathbb{F}_\delta$; it
can therefore be considered as the pullback via $\eta$ of a 
divisor on $\mathbb{F}_\delta$
of the form 
\begin{equation}\label{graph of bundle}
\mathcal{G}_E := \sum_{i=1}^k f_i + A,
\end{equation}
where $f_i$ is the fibre of the ruled surface $\mathbb{F}_\delta$ over 
the point $x_i$ and
$A$ is a section of the ruling such that $\eta^\ast A = \overline{C}$. 
The divisor $\mathcal{G}_E$ is called the {\em graph} of the bundle $E$.
We finish by noting that, 
although the section $A$ is a smooth curve on $\mathbb{F}_\delta$,
its pullback need not be smooth: it may be reducible or multiple with
multiplicity 2.
\begin{remark*} 
If $\delta$ is the pullback of a line bundle on $B$, then 
its restriction to any fibre 
of $\pi$ is trivial and the induced involution $i_\delta$ 
is given by $(b,\lambda) \mapsto (b,\lambda^{-1})$; 
in this case, we have $(B \times T^\ast)/i_\delta = B \times \mathbb{P}^1$.
Furthermore, if there exist line bundles $a$ and $\delta'$ on $X$ such that 
$\delta = a^2\delta'$, then 
$\mathbb{F}_\delta$ is isomorphic to $\mathbb{F}_{\delta'}$;
indeed, the map $a: J(X) \rightarrow J(X)$ defined by 
$(b,\lambda) \mapsto (b,a_b \lambda)$ is an isomorphism of the Jacobian surface 
that commutes with the involutions determined by $\delta$ and $\delta'$. 
In particular, if $\delta$ is an element of $2NS(X)$, then $\delta = a^2$ for some 
line bundle on $X$ and $\mathbb{F}_\delta$ is isomorphic to $B \times \mathbb{P}^1$. 
\end{remark*}

For any $c_1$ in $NS(X)$, choose a line bundle $\delta$ on 
$X$ such that 
\[ c_1(\delta) \in c_1 + 2NS(X)\] 
and 
\[ m_{c_1} := m(2,c_1) = -\frac{1}{2} \left( c_1(\delta )/2 \right)^2. \] 
Therefore, if $\delta'$ is any other line bundle with Chern class in 
$c_1 + 2NS(X)$, it induces a ruled surface that is isomorphic to 
$\mathbb{F}_\delta$; the advantage of using this particular $\delta$ is that its 
Chern class has maximal self-intersection $-8m_{c_1}$.

Let us now compute the invariant of the ruled surface.
We begin by setting some notation.
We denote $B_0$ the zero-section of $J(X)$ and 
$\Sigma_\delta$ the section in $J(X)$ corresponding to $\delta$; 
also, let $p_1: J(X) \rightarrow B$ be
the projection onto the first factor. Consider the exact sequence
\[ 0 \rightarrow \mathcal{O}_{J(X)}(\Sigma_\delta) \rightarrow 
\mathcal{O}_{J(X)}(B_0 + \Sigma_\delta) \rightarrow \mathcal{O}_{B_0}
(\Sigma_\delta) \rightarrow 0. \]
Pushing down to $B$, we obtain a new exact sequence
\begin{equation}\label{exact1}
0 \rightarrow \mathcal{O}_B \rightarrow V_\delta \rightarrow L \rightarrow 0,
\end{equation}
where 
\[ V_\delta := {p_1}_\ast(\mathcal{O}_{J(X)}(B_0 + \Sigma_\delta))\]
is a rank-2 vector bundle on the curve $B$ and 
\[ L := {p_1}_\ast(\mathcal{O}_{B_0}(\Sigma_\delta))\] 
is a line bundle of degree $4m_{c_1}$ on $B$, given by the effective divisor 
$q_1 + \dots + q_{4m_{c_1}}$ that corresponds to the projection onto 
$B$ of the intersection points $B_0 \cap \Sigma_\delta$ (counted with multiplicity). 
Note that $\mathbb{F}_\delta = \mathbb{P}(V_\delta)$.
Given the above notation, we have the following result. 
 
\begin{lemma}\label{ruled surface invariant} 
Let $d$ be the maximal degree of a
subline bundle of $V_\delta$; it is then a non-negative integer that satisfies the inequality
\[ \max \{ 0 , 2m_{c_1} - g/2 \} \leq d \leq 2m_{c_1}. \]
Moreover, the invariant of the ruled surface 
$\mathbb{F}_\delta = \mathbb{P}(V_\delta)$ is
\[ e = 2d - 4m_{c_1}. \] 
\end{lemma}

\begin{proof}
The invariant $e$ of the ruled 
surface is given by
\[e = \max \{ \mbox{$2\deg{\lambda} - \deg{V_\delta}$ : 
there exists a nonzero map $\lambda \rightarrow V_\delta$} \}, \]
where $\lambda$ is a line bundle on $B$ (see, for example, \cite{F2}). 
Therefore, if $d$ is the maximal degree of a subline bundle of $V_\delta$,
we have 
\[ e = 2d - \deg{V_\delta} = 2d - 4m_{c_1}. \]
Note that, since $\mathcal{O}$ is a subline bundle of $V_\delta$ (see \eqref{exact1}),
the integer $d$ is non-negative.
To determine the bounds of $d$, we have to verify that 
\[ -g \leq e \leq 0.\] 
The left-hand inequality follows from a theorem of Segre-Nagata \cite{F2};
hence, there only remains to show that $e$ is less than or equal to zero. 

Let $A$ be a section of the ruled surface $\mathbb F_\delta$;
the pullback $\eta^*A$ is therefore a bisection of $J(X)$. If it is 
reducible, then its two components are sections $C_1$ and $C_2$ of $J(X)$, 
giving
\[2A^2=(\eta^*A)^2=(C_1+C_2)^2=2C_1\cdot C_2\geq 0.\]
If the bisection ${\overline C} = \eta^*A$ is instead irreducible, we consider
its normalization $C\rightarrow {\overline C}$ and let $\gamma :C\rightarrow B$ 
be the two-to-one map induced by $C\rightarrow {\overline C}
\subset J(X)$.
Note that the natural map $C \rightarrow J(X) \times_B C$ gives a section $C_1$ of the surface 
$C \times T^\ast \rightarrow C$;
moreover, if we denote by ${\tilde\gamma}:C\times T^* \rightarrow J(X)$
the two-to-one map induced by $\gamma$, then the pullback
${\tilde\gamma}^*({\overline C})$ is reducible, with components $C_1$ and $C_2$,
where $C_2$ is also section of $C\times T^*\rightarrow C$, and we have
\[ 4A^2= ({\tilde\gamma}^*({\overline C}))^2 = (C_1+C_2)^2 = 2C_1 \cdot C_2\geq 0.\]
Therefore, since 
\[ e=-\min \{ \mbox{$A^2$ $|$ $A$ section of $\mathbb{F}_\delta$} \}\]
(see \cite{F2}, Proposition 12, Chapter 5), 
it follows that $e$ is non-positive. 
\end{proof}

\begin{remark*}
For a generic curve $B$ of genus greater than $1$, 
the Neron-Severi group of an elliptic surface $X$ over $B$ is trivial 
and the ruled surface is $B \times \mathbb{P}^1$ for any $\delta$ in $\text{Pic}(X)$. 
Moreover, this is always true if $B$ is rational:
the sections of the ruled surface are given by rational maps 
$\mathbb{P}^1 \rightarrow \mathbb{P}^1$ and the irreducible bisections of $J(X)$
are the pullbacks to $J(X)$ of non-constant rational maps (for details, see \cite{Moraru}).
\end{remark*}

We finish this section by determining the genus of irreducible bisections.
\begin{lemma}
If the spectral cover of the bundle $E$ is a smooth irreducible bisection 
$\overline{C}$ of $J(X)$, then its genus is given by
\begin{equation}\label{genus of spectral curve}
g(\overline{C}) = 4\Delta(E) + 2g - 1, 
\end{equation}
where $g$ is the genus of $B$.
\end{lemma}

\begin{proof}
We begin by noting that the pushforward $A_0 := \eta_\ast(B_0)$ of the zero section of $J(X)$
is a section of the ruled surface $\mathbb{F}_\delta$ whose
pullback $\eta^\ast A_0$ to $J(X)$ is the reducible bisection $B_0 + \Sigma_\delta$;
consequently, it has self-intersection $A_0^2 = - c_1^2/2$.
We now describe the ramification and branching divisors of $\eta$. 
Let $R$ be the ramification divisor in $J(X)$, defined as
the fixed point set of $\eta$;
referring to Lemma \ref{degree of section}, we have
\[ R \cdot B_0 = \# \{ (b,t) \ : \ \delta_b  = \mathcal{O}_T \} = \Sigma_\delta \cdot B_0 = - c_1^2/2. \]
The branching divisor $G$ is a 4-section of $\mathbb{F}_\delta$
such that $\eta^\ast G = 2R$;
since
\[ G \cdot A_0 = G \cdot \eta_\ast(B_0) = \eta_\ast(\eta^\ast G \cdot B_0) = -c_1^2, \]
it is equivalent to a divisor of the form $4A_0 + \mathfrak{b} f$,
where $\mathfrak{b}$ is a divisor on $B$ of degree $c_1^2$
and $f$ is a fibre of the ruled surface. 

Let $A$ be the graph of the bundle $E$, that is, the section of $\mathbb{F}_\delta$ such that  
$\overline{C} = \eta^\ast A$.
If we write $A \sim A_0 + \mathfrak{b}'f$, for some divisor $\mathfrak{b}'$ on $B$, then
$\overline{C} \sim (B_0 + \Sigma_\delta) + \mathfrak{b}' T^\ast$, 
where $\mathfrak{b}'$ also denotes the pullback of the divisor to $J(X)$,
and the intersection number $\overline{C} \cdot B_0$ is equal to $-c_1^2/2 + \deg{\mathfrak{b}'}$.
Recall that $\overline{C} \cdot B_0$ is, by construction, the number of points 
(counted with multiplicity) in the support of the torsion sheaf 
$R^1\pi_\ast(E \otimes \mathcal{O}_X)$, which is equal to 
$n_E = c_2 - c_1^2/2$ (see section \ref{spectral curve}); therefore, we have $\deg{\mathfrak{b}'} = c_2$.
Hence, the smooth bisection $\overline{C}$ is a double cover of $B$ of branching order 
$G \cdot A = 4c_2 - c_1^2$ and \eqref{genus of spectral curve} follows by the Hurwitz formula. 
\end{proof}

\section{Existence theorems} 
Let $E$ be a holomorphic rank-2 vector bundle on the non-K\"{a}hler elliptic surface
$X$ with determinant line bundle $\delta$ and Chern classes $c_1$ and $c_2$. 
If we denote
\[ \Delta (E) := \Delta(2,c_1,c_2)\]
the discriminant of $E$, then a well-known result states 
that $\Delta (E)$ cannot be negative \cite{BaL,ElFo,BrF,B3,LeP}.

\subsection{Rank-2 vector bundle as extensions}
\label{existence theorems}

By using Lemma \ref{degree of restriction of line bundle}, 
Proposition \ref{section of line bundle} and 
Lemma \ref{degree of section}, one obtains the following 
result, whose proof is similar to that of Theorem 1.3, Chapter VII, \cite{FM}: 
\begin{theorem}\label{bound}
Let $\pi :X\rightarrow B$ be a non-K\"{a}hler elliptic
surface and $E$ be a holomorphic rank-2 vector bundle on $X$ 
with determinant line bundle $\delta$. 
Then $E$ satisfies one of the following two cases:

(A) There exists a line bundle $\mathcal{D}$ on $X$ and a locally complete
intersection $Z$ of codimension 2 in $X$ such that $E$ is given by an
extension 
\[0\rightarrow \mathcal{D}\rightarrow E\rightarrow \delta\otimes 
\mathcal{D}^{-1} \otimes I_Z\rightarrow 0.\]
In fact, $Z$ is the set of points (counted with multiplicity) corresponding to the 
fibres of $\pi$ over which the bundle $E$ is unstable.
Moreover, we have 
\[ \Delta (E)=\frac{1}{8}{\overline C}^2+\frac{1}{2}\ell (Z).\]

\indent
(B) There exists:
(i) a smooth irreducible curve $C$ and a birational map $C\rightarrow 
\overline{C}\subset J(X)$, where $\overline C$ is a bisection
that is invariant under the involution $i_\delta$ on $J(X)$ defined by 
the line bundle $\delta$;

(ii) a line bundle $\tilde{\mathcal{D}}$ on the normalisation $W$ of 
$X\times_BC$, whose restriction to a smooth fibre of 
$W\rightarrow C$ is the same as the one induced by the section of $J(W)$ that
corresponds to the map $C\rightarrow J(X)$;

(iii) a codimension 2 locally complete intersection $\tilde Z$ in $W$,
an exact sequence
\[0\rightarrow \tilde{\mathcal{D}}\rightarrow \tilde{\gamma}^*E\rightarrow
\tilde{\gamma}^*\delta \otimes \tilde{\mathcal{D}}^{-1}\otimes I_{\tilde Z}
\rightarrow 0,\]
where $\tilde{\gamma} :W\rightarrow X$ is the natural map, 
and 
\[ \Delta (E)=\frac{1}{8}{\overline C}^2+\frac{1}{4}\ell ({\tilde Z}). \]
This time, $\tilde{Z}$ is the set of points corresponding to the 
fibres of $W \rightarrow C$ over which the bundle ${\tilde\gamma}^*E$ is unstable.
\hfill $\square$
\end{theorem}

\begin{remark}\label{reducible bundles}
Suppose that the vector bundle $E$ satisfies case (A) of Theorem \ref{bound}. 
Let $\Sigma_1$ and $\Sigma_2$ be the
sections of $J(X)$ determined by the line bundles $\mathcal{D}$ and 
$\mathcal{D} \otimes \delta$, respectively.
Then, one can easily verify that $\overline{C} = \Sigma_1 + \Sigma_2$, implying that the
bisection associated to $E$ is reducible or a section counted with multiplicity 2 
(if $\Sigma_1 = \Sigma_2$).
\end{remark}

We now have the following complete description of  non-filtrable bundles: 

\begin{proposition} 
Let $E$ be any holomorphic $2$-vector bundle over $X$. 
Suppose that the spectral cover of $E$ includes the bisection 
$\overline{C}$ of $J(X)$. Then $E$ is non-filtrable if and only if $\overline{C}$ is irreducible.
\end{proposition}

\begin{proof}
Suppose that there exits a line bundle $\mathcal{D}$ on $X$ that maps into 
$E$. After possibly tensoring $\mathcal{D}$ by the pullback of a suitable 
line bundle on $B$, the rank-2 bundle $E$ is then given as an extension
\[ 0 \rightarrow \mathcal{D} \rightarrow E \rightarrow 
\mathcal{D}^{-1} \otimes \delta \otimes I_Z \rightarrow 0, \]
where $Z \subset X$ is a locally complete intersection of codimension 2,
that is, $E$ satisfies case (A) of Theorem \ref{bound}; 
referring to remark \ref{reducible bundles}, the bisection is then not irreducible. 
Conversely, suppose that the bisection is not irreducible and that $\Sigma$ is one of its components.
If $\mathcal{D}$ is a line bundle on $X$ corresponding to $\Sigma$, then $\mathcal{D}$ 
maps non-trivially into $E$, implying that $E$ is filtrable.
\end{proof}

\begin{note}
A partial characterisation of non-filtrable bundles is also given in \cite{AT}.
\end{note}

\subsection{Existence of rank-2 vector bundles}
\label{existence theorems: general case}

A partial converse of Theorem \ref{bound} is the following result:
\begin{theorem}\label{existence}
Let $\pi :X\rightarrow B$ be a non-K\"{a}hler elliptic surface and 
$\delta$ be a line bundle in $Pic(X)$. 
Furthermore, let $i_\delta:J(X)\rightarrow J(X)$ 
be the involution defined by $\delta$ and suppose that 
$\overline C$ is a bisection of $J(X)\rightarrow B$
that is invariant with respect to the involution $i_\delta$. Then, there exists 
a rank-2 holomorphic vector bundle $E$ on $X$ such that 
\[\mbox{$c_1(E)=c_1(\delta )$ and $\Delta (E)=\dfrac{1}{8}{\overline C}^2=
\dfrac{1}{4}A^2$},\]
where $A$ is a section of the ruled surface $\mathbb F_\delta$ with $\eta^*A=\overline C$.
\end{theorem}

\begin{proof} 
If the bisection $\overline{C}$ is reducible, then its components are 
sections $\Sigma_1$ and $\Sigma_2$ of $J(X)$. Let $\mathcal{D}$ be a line bundle
on $X$ corresponding to $\Sigma_1$ (see Proposition
\ref{section of line bundle}); if $E$ is any extension of 
$\mathcal{D}^{-1} \otimes \delta$ by $\mathcal{D}$, then $E$
is a rank-2 vector bundle on $X$ that has determinant $\delta$ and spectral cover $\overline{C}$.

If the bisection $\overline{C}$ is irreducible,
then consider its normalisation $C\rightarrow {\overline C}$ 
and let $\gamma :C\rightarrow B$ be the double covering induced 
by $C\rightarrow {\overline C}\subset J(X)$. The
normalisation $W$ of the fibred product $X\times_B C$
is then a non-K\"{a}hler elliptic surface over $C$ with relative
Jacobian $J(W)=C\times T^*$; furthermore, the natural two-to-one map 
$\tilde\gamma :W\rightarrow X$ induces a covering 
$\gamma ':J(W)\rightarrow J(X)$. 
Note that the inclusion map $C\rightarrow J(X)\times_B C$ gives a 
section $\Sigma_1$ of $J(W)\rightarrow C$; 
the pullback ${\gamma '}^*{\overline C}$
is then reducible with components $\Sigma_1$ and $\Sigma_2$, 
where $\Sigma_2$ is another section of $J(W)$. 
By Proposition \ref{section of line bundle}, there exists 
a line bundle $\mathcal{L}$ on $W$ whose restriction to any smooth 
fibre $T_c$ of $W$ is $\Sigma_{2c}$.
Let $\mathcal{D}$ be the line bundle on $W$
satisfying the equality
\[\mathcal{L}\cong {\tilde\gamma}^*\delta \otimes \mathcal{D}^{-1}\]
and define the holomorphic rank-2 vector bundle $E$ on $X$ by
\[E:={\tilde\gamma}_*(\mathcal{ L});\]
we then have to show that $E$ has first Chern class $c_1(\delta)$ and 
discriminant $\frac{1}{8}\overline{C}^2$. 

Let $\overline{i}_\delta$ be the involution on $W$ 
that interchanges 
the sheets of $\tilde\gamma$. If $G \subset X$ is the (smooth) 
branch divisor of the 
double covering $\tilde\gamma :W\rightarrow X$, then there exists a line 
bundle $L_0$ on $X$ such that $L_0^2=\mathcal{O}_X(G)$; 
moreover, by Lemma 29, Chapter 2 of \cite{F2} or by \cite{B4}, there is an 
exact sequence:
\[0\rightarrow {\overline{i}_\delta}^*\mathcal{L}\otimes {\tilde\gamma}^*L_0^{-1}
\rightarrow {\tilde\gamma}^*{\tilde\gamma}_*(\mathcal{L})\rightarrow \mathcal{L}
\rightarrow 0.\]
Since the involution $\overline{i}_\delta$ on $W$ is induced by interchanging 
the sheets of the double cover $C\rightarrow B$,  
the restriction of $\overline{i}_\delta^*\mathcal{L}$ 
to any smooth 
fibre $T_c$ of $W$ (which is not in the ramification locus of $\tilde\gamma$)
is isomorphic to the restriction of $\mathcal{D}$ to the same fibre, namely to $\Sigma_{1c}$. 
From the preceding exact sequence, we obtain
\[0\rightarrow \mathcal{D}\otimes\mathcal{O}_W(F)\rightarrow {\tilde\gamma}^*E
\rightarrow {\tilde\gamma}^*\delta\otimes \mathcal{D}^{-1}\rightarrow 0,\]
where $F$ is a divisor on $W$ (hence a combination of fibres of the 
non-K\"{a}hler elliptic surface $W\rightarrow C$). 
Referring to Theorem \ref{bound}, we have 
\[ \Delta (E)=\frac{1}{8}\overline{C}^2 = \frac{1}{4}A^2,\]
where $A$ is the section of the ruled surface $\mathbb F_\delta$ defined by 
the bisection $\overline C$. By \cite{ABT}, we also have
\[c_1(E)\equiv c_1(\delta )\;\;\mbox{mod}\;\mbox{Tors}\;(NS(X)).\]
To get rid of the torsion, we need to add multiples of classes of fibres. 
Then, as in \cite{ABT}, we can modify the line bundle $\mathcal{L}$, by tensoring 
it with line bundles of the form $\mathcal{O}_W(T_c)$ or $\mathcal{O}_W(T_i)$, 
and obtain the desired result
\[c_1(E)=c_1(\delta ).\]
Note that the discriminant remains unchanged (see the formula in \cite{ABT} 
for the direct image of a line bundle).
\end{proof}

The above result implies that the existence problem for vector bundles is equivalent to 
the existence problem of bisections of $J(X)$ that are invariant under a given involution. 
Let us fix an element $c_1$ in $NS(X)$ and a line bundle $\delta$ on $X$ such 
that $c_1(\delta) \in c_1 + 2NS(X)$
and $m_{c_1} := m(2,c_1) 
= -\frac{1}{2} \left( c_1(\delta )/2 \right)^2$. 
Referring to section \ref{graph} and Lemma \ref{ruled surface invariant},
the Jacobian surface $J(X)$ of $X$ is 
thus endowed with an involution $i_\delta$ and the quotient 
is a ruled surface $\mathbb{F}_\delta$
that has a non-positive invariant $e$; 
moreover, there is a one-to-one correspondence between sections of $\mathbb{F}_\delta$ and 
spectral curves of rank-2 vector bundles on $X$ that have determinant $\delta$ and no jumps. 
Therefore, the minimum value of the discriminant of a vector bundle 
$E$ on $X$ with first Chern class $c_1$ is equal to $-e/4$. 
Conversely, one can show that for any integer $c_2$ such that  
$\Delta(2,c_1,c_2)$ is greater or equal to $-e/4$,
there exists a rank-2 vector bundle on $X$ with Chern classes $c_1$ and $c_2$. 
We can now state the main result of the paper:
\begin{theorem}\label{general existence}
Let $X$ be a minimal non-K\"{a}hler elliptic surface over a curve $B$ of genus 
$g$ and fix a pair $(c_1,c_2)$ in $NS(X) \times \mathbb{Z}$. 
Let $m_{c_1} := m(2,c_1)$ and choose a line bundle $\delta$ on $X$ such that 
$-c_1^2(\delta)/2$ is equal to $4m_{c_1}$. 
Then, there exists a holomorphic rank-2 vector bundle on $X$ with 
Chern classes $c_1$ and $c_2$ if and only if 
\[\Delta(2,c_1,c_2) \geq (m_{c_1}-d/2),  \]
where $d$ is the non-negative integer determined in Lemma \ref{ruled surface invariant}. 
Furthermore, if 
\[ (m_{c_1}-d/2) \leq \Delta(2,c_1,c_2) < m_{c_1}, \]
then the corresponding vector bundles are non-filtrable.
\end{theorem}

\begin{proof}
Recall from Lemma \ref{ruled surface invariant} that the invariant $e$ of the 
ruled surface is equal to $2d - 4m_{c_1}$. Let $\Delta_0 := -e/4 = m_{c_1} - d/2$
and consider $\Delta := \Delta(2,c_1,c_2)\geq \Delta_0$; 
note that $k=2(\Delta -\Delta_0) \geq 0$ is an integer.
It is sufficient to prove the existence of a holomorphic rank-$2$ vector 
bundle $E$ with first Chern class $c_1(\delta)$ and discriminant $\Delta$. 
Let $\overline{C}_0$ be a bisection of $J(X)$ of minimal self-intersection $8\Delta_0$.
If $k=0$, choose a holomorphic rank-$2$ vector bundle $E_0$ corresponding to $\overline{C}_0$,
for example, any bundle determined by Theorem \ref{existence}.

For $k > 0$, choose a smooth fibre $T := \pi^{-1}(b)$ of $\pi$, with $b \in B$,
such that if the bisection $\overline{C}_0$ is irreducible, 
then the double cover $\overline{C}_0 \rightarrow B$ does not have a branch point over $b$.
Set $\delta ':= \delta \otimes {\mathcal O}_X(kT)$.
The line bundles $\delta$ and $\delta'$ then both correspond to the same section in $J(X)$,
inducing isomorphic ruled surfaces ${\mathbb F}_{\delta '}$ and ${\mathbb F}_{\delta}$, respectively.
Consequently, there exists a holomorphic rank-$2$ 
vector bundle $E'_0$ on $X$ with first Chern class $c_1(\delta ')$ and 
discriminant $\Delta_0$ that is regular on the fibre $T$
(over an elliptic curve, a bundle is said to be {\em regular} if
its group of automorphisms is of the smallest possible dimension). 
Indeed, if $\overline{C}_0$ is reducible, then choose line bundles $L_1$ and $L_2$ on $X$ 
associated to the components of $\overline{C}_0$, with $L_1 \otimes L_2 = \delta'$, and let
$E'_0$ be an extension of $L_2$ by $L_1$ that is regular on $T$.
Moreover, if $\overline{C}_0$ is irreducible, then $E'_0$ can be any vector
bundle given by Theorem \ref{existence}.  
Let $j:T\rightarrow X$ be the natural inclusion map;
if $\lambda$ is a line bundle on $T$ of degree $1$, then there exists
a surjection $E'_0\rightarrow j_*\lambda$.
Consider the elementary modification
\[ 0\rightarrow E_1\rightarrow E'_0\rightarrow j_*\lambda\rightarrow 0; \]
then, the bundle $E_1$ splits as $\lambda \oplus \lambda^\ast$ over $T$ and 
there exists a surjection 
$E_1 \rightarrow j_*\lambda$. Hence, by performing $(k-1)$ successive elementary 
modifications on $E_1$ with respect to $j_\ast\lambda$, one obtains a holomorphic 
vector bundle $E$ on $X$ with first Chern class $c_1(\delta )$ and 
discriminant $\Delta$. 
\end{proof} 

\begin{remark*}
If the genus of the base curve $B$ is less than $2$, then the statement of the theorem becomes: 
there exists a holomorphic rank-2 vector bundle $E$ on $X$ with
Chern classes $c_1$ and $c_2$ if and only if the discriminant $\Delta(2,c_1,c_2)$
is a non-negative number. 
(For an alternate proof in the case of primary Kodaira surfaces, see \cite{ABT}.)
In contrast, if the genus of the base curve is greater  $1$, there are "gaps" for the discriminant 
of holomorphic rank-$2$ vector bundles, whenever $m_{c_1}$ is greater than $d/2$; thus, 
the existence of holomorphic 
vector bundles on $X$ depends on the geometry of the base curve $B$. 
However, by the proof of Theorem \ref{general existence},
once there is an irreducible bisection of $J(X)$, one can construct 
infinitely many non-filtrable vector bundles.
\end{remark*}
\begin{note}
Bundles with $\Delta (2,c_1,c_2)=0$ have also
been studied in \cite{A-B}. 
\end{note}

{\bf Acknowledgements}
The first author would like to express his gratitude 
to the Max Planck Institute of Mathematics for its hospitality and 
stimulating atmosphere; this paper was prepared during his stay at the 
Institute. It is a pleasure for both authors to thank Jacques Hurtubise for 
suggesting a link between the papers \cite{ABT} and \cite{Moraru}.
The second author would like to thank Jacques Hurtubise for his generous 
encouragement and support during the completion of this paper.
She would also like to thank Ron Donagi and Tony Pantev for valuable 
discussions, and the
Department of Mathematics at the University of Pennsylvania for their 
hospitality, during the preparation of part of this article.

\end{document}